\documentstyle[12pt]{article}
\newcommand{\Reals}{\rm I\kern-.19emR}
\newcommand{\Notin}{/\kern-.6em\hbox{$\in$}}
\newcommand{\MM}{\rm I\kern-.19emM}
\newcommand{\Notequiv}{/\kern-.6em\hbox{$\equiv$}}
\newcommand{\Ceals}{\rm I\kern-.5emC}
\newcommand{\nsubset}{/\kern-.6em\hbox{$\subset$}}

\newcommand{\nin}{\backslash \kern-.5em\in}
\newtheorem{theo}{Theorem}
\newtheorem{coro}{Corollary}
\newtheorem{prop}{Proposition}
\newtheorem{defn}{Definition}
\begin{document}
\begin{center}
{{\large \bf ON METHOD OF STATISTICAL}\\[1ex]
{{\large \bf DIFFERENTIALS}}\\[1ex]
{by}\\[1ex]
{Rohitha Goonatilake\footnote{
 Current Address:
 Department of Mathematical and Physical Sciences, 
 Texas A\&M International University, Laredo, Texas 78041-1900.\\
 E-mail: harag@tamiu.edu}}}
\end{center}
\vskip 1em
\noindent{\bf ABSTRACT:~}
 The method of statistical differentials, which approximates the mean value and
 variance of transformations of random variables is used in many areas of
 mathematics. This paper will discuss the conditions under which such an 
 approximation will be exact, and also explore their accuracy in terms of 
 error bounds under certain moment conditions.
\vskip 1em
\noindent{{\bf AMS SUBJECT CLASSIFICATION (2000):}~ 40A05, 41A10, 41A25, 60A99}\\
\noindent{{\bf KEYWORDS:}~ Approximation, Error, Exactness, Expected value, and Variance.}\\
\vskip 1em
\noindent{\bf I.~ PRELIMINARIES:}
\vskip 1em

 The method of statistical differentials is a method of approximating
 the mean (expected value), and the variance of transformations of the
 random variables. The method has been presented by Johnson and Kotz [1],
 Elandt-Johnson and Johnson [2], and London [3]. Frye [4] has given 
 counter examples for the approximation to be exact, and also to show that
 in some cases, this approximation will not hold beyond certain degreed
 polynomials.

 Throughout this paper, we will use the notations adopted by London [3].
 Let $Y = g(X),$ where $X$ is a random variable such that $E[X]$ and 
 $\hbox{Var}[X]$ exist and are known. Other moment conditions will be assumed, 
 depending upon the cases, that are discussed from time to time. Let
 $m = E[X].$ In order to express $Y = g(X)$ as a Taylor series expanded
 about $X = m,$ assume $g(X)$ is a function possessing derivatives of all 
 order up to $n^{\hbox{th}}$ throughout the interval $a \leq X \leq b,$ then
 there is a value $z,$ with $z \in (m, X) \cup (X, m)$ such that 
\begin{eqnarray*}
   g(X) & = & g(m) + (X -m)g^\prime(m) + 
              {(X -m)^2\over 2!}g^{\prime\prime}(m) + \cdots\\
 \cdots & + & {(X -m)^{n -1}\over (n -1)!}g^{(n -1)}(m) + 
            {(X -m)^n\over n!}g^{(n)}(z).
\end{eqnarray*}
 The random variable $z$ lies between $m$ and $X.$ Taking expected values
 of both sides gives,
\begin{eqnarray*}
 E[Y] & = & E[g(X)] = g(m) + E[X -m]g^\prime(m) + 
          E\Biggl[{(X -m)^2\over 2!}\Biggr]g^{\prime\prime}(m) + \cdots\\
 \cdots & + & E\Biggl[{(X -m)^{n -1}\over (n -1)!}\Biggr]g^{(n -1)}(m) + 
            E\Biggl[{(X -m)^n\over n!}g^{(n)}(z)\Biggr].
\end{eqnarray*}
 Since $E[X -m] = 0$ and $E[(X -m)^2] = \hbox{Var}[X],$ we have
\begin{eqnarray}\label{1}
 E[Y] & = & g(m) +  {g^{\prime\prime}(m)\over 2!}\cdot\hbox{Var}[X] +
          {E[(X -m)^3]\over 3!}g^{\prime\prime\prime}(m) + \cdots\nonumber\\
 \cdots & + & E\Biggl[{(X -m)^{n -1}\over (n -1)!}\Biggr]g^{(n -1)}(m) + 
            E\Biggl[{(X -m)^n\over n!}g^{(n)}(z)\Biggr],
\end{eqnarray}
 where $n$ is a positive integer greater than 1. It is customary to truncate 
 this series at the second term, and to consider
\begin{eqnarray*}
 E[Y] \cong g(m) +  {g^{\prime\prime}(m)\over 2!}\cdot\hbox{Var}[X]
\end{eqnarray*}
 as statistical differential approximation. We refer to the
 approximation (\ref{1}) as a statistical differential approximation.

\vskip 1em
\noindent{\bf II.~MOTIVATION:}
\vskip 1em

 We first give an important order notation, which we use time to time 
 throughout this paper. The $O-$notation (read {\sl big-oh notation}) 
 provides a special way to compare relative sizes of functions that is
 very useful in the analysis of error bounds. The $o-$notation (read {\sl
 small-oh notation} is given for completeness of the definition and
 will not be used in the discussion elsewhere.

\begin{defn}{\rm
 {\bf Landau Order Notations:}\\
 Let $f(x)$ and $g(x)$ be given functions. Let $x_0$ be a fixed point and
 suppose that $g(x)$ is positive and continuous in an open interval about
 $x_0,$ where $x_0$ may be finite or infinite.
\begin{enumerate}
 \item If there is a constant $K$ such that 
       $$|f(x)| < Kg(x)$$
       in an open interval about $x_0,$ then $f(x) = O\Bigl(g(x)\Bigr),$
       $(x \to x_o).$
 \item Furthermore, if 
       $$\lim_{x \to x_0}{f(x)\over g(x)} = 0,$$
       then $f(x) = o\Bigl(g(x)\Bigr),$
       $(x \to x_o).$
\end{enumerate}}
\end{defn}

 Three conditions (i), (ii) and (iii) below, each of them will essentially 
 lead to exactness of the approximation. It is reminded that (i) and
 (ii) can not be true. For (ii), a random variable with a symmetric distribution
 around 0, $m = 0$ and $EX^j = E[(X - m)^j] = 0$ for all odd $j$ (and in
 particular $j = 3,$) but not for even $j.$

\begin{description}
 \item (i). 
       $g^{\prime\prime\prime}(m) = 0$ implies that all other derivatives
       of order greater than three evaluated at $x = m$ to be zero.
       In particular, $g^{(n)}(z) = 0$ for the random variable $z$ lies 
       between $m$ and $X.$
 \item (ii).
       $E\bigl[(X -m)^3\bigr] = 0,$ implies that all other central moments
       of order greater than three about mean equal zero.
 \item (iii).
       Remaining terms beyond third sum up to a zero.
\end{description}

 We use the following polynomial expansion later in the paper. The series expansion 
 $(x + y)^n$ is symmetric with respect to the changes of variables $x, y,$ so does 
 convergence region. The expression in parenthesis following of the series, indicates
 the region of convergence. If not otherwise indicated, it is to be understood that 
 the series converges for all values of the variable.
 $$(x + y)^n = x^n + n x^{n -1}y + {n(n -1)\over 2!}x^{n -2}y^2 +
               {n(n -1)(n -2)\over 3!}x^{n -3}y^3 + \cdots, ~\hbox{where}~
                y^2 < x^2.$$

\vskip 1em
\noindent{\bf III.~EXACTNESS OF THE EXPECTED VALUE:}
\vskip 1em

 In this section, attention is drawn to the conditions for which an 
 exactness of the approximation can be achieved.

\begin{prop}
 If $g(\cdot)$ is a polynomial of degree 2, then the statistical
 differential approximation for the expected value of the transformations
 of the random variable $X,$ $E\bigl[g(X)\bigr]$ is exact.
\end{prop}

 The above proposition is extended to involve third degree polynomial
 in the following manner.

\begin{prop} 
 If $g(\cdot)$ is a polynomial of degree 3, then the statistical
 differential approximation for the expected value of the transformations
 of the random variable $X,$ $E\bigl[g(X)\bigr]$ is exact, provided 
 $E\bigl[(X -m)^3\bigr] = 0.$
\end{prop}

 The natural question is that could this procedure be so extended to the
 next higher degree polynomial by requiring vanishing fourth central moment
 about mean of the random variable, (in addition to what have been already 
 assumed). A more generalized version of the above proposition is the
 following.

 Suppose $g(\cdot)$ is a polynomial of degree $j,$ then requiring
 $E\bigl[(X -m)^k\bigr] = 0,$ for all $3 \leq k \leq j$ will do the job! 
 But, for $k = 4,$ this gives that $E\bigl[(X -m)^4\bigr] = 0,$ which
 immediately implies that we are dealing with a constant random variable, and
 there is no need to worry about anything else.

 The condition of the last assertion leads to a nice relation of the 
 $3\hbox{rd}$ moment of the random variable $X.$

\begin{theo}
 If $X$ is a random variable such that its mean and variance exist and are 
 known, together with the property,
 $E\bigl[(X -m)^3\bigr] = 0,$ then $3\hbox{rd}$
 moment of the random variable $X,$ $E[X^3]$ exists, is finite, and satisfies 
 the relation,
 $$E[X^3] = \bigl(E[X]\bigr)^3\Biggl\{3\cdot{\hbox{Var}[X]\over
            \bigl(E[X]\bigr)^2} + 1 \Biggr\}.$$
\end{theo}

\noindent{\bf Proof.} 

 From the series expansion for $(x + y)^n$ with $y^2 < x^2,$ and letting
 $x = X - m,$ and $y = m > 0,$ where $E[X] =m,$ we obtain,
 $$X^j = (X - m)^j + j(X -m)^{j -1}m + {j(j -1)\over 2!}(X - m)^{j -2}m^2 +
          {j(j -1)(j -2)\over 3!}(X - m)^{j -3}m^3 + \cdots$$
 $$     + {j(j -1)(j -2)\cdots 3\over (j -2)!}(X - m)^2m^{j -2} +
          {j(j -1)(j -2)\cdots 2\over (j -1)!}(X - m)m^{j -1} + m^j.$$
 For $j =3,$ taking the expected value, and using  
 $E\bigl[(X -m)^3\bigr] = 0,$ we obtain,
 $$E[X^3] = \bigl(E[X]\bigr)^3\Biggl\{3\cdot{\hbox{Var}[X] \over
            \bigl(E[X]\bigr)^2} + 1 \Biggr\},$$
 as asserted.

\vrule height 8pt width 4pt

 Our next task would be to see, are there other conditions for which this 
 approximation is exact ?
 One of the results in this connection is to consider, the {\sl Peano Kernel}
 method [5];

 For any $g \in C^{n +1}[a, b],$ the Taylor expansion with integral remainder
 gives,
 $$g(X) = \sum_{k = 0}^n {(X - a)^k\over k!}g^{(k)}(a) + 
          {1\over n!}\int_a^X (X - \theta)^ng^{(n +1)}(\theta)d\theta,
          ~\hbox{where}~a \leq X \leq b,$$ and $g \in C^{n +1}[a, b]$ means 
 $g(\cdot)$ is $(n +1)$ continuously differentiable function over $[a, b].$
 Based on the above formula, we have:

\begin{theo}
 If $g(\cdot)$ has derivatives of order $j \leq 3,$ $g \in C^3[I],$ 
 where $I$ is some interval containing the range of $X,$ and
 $E\Bigl(\int_m^X (X - \theta)^2 g^{\prime\prime\prime}(\theta)d\theta\Bigr) = 0,$ 
 then the statistical differential approximations for the expected value 
 of the transformations of the random variable $X,$ $E[g(X)]$ is exact.
 Otherwise, it will be exact up to the error term
 $O\Bigl(E\Bigl(\int_m^X(X - \theta)^2
         g^{\prime\prime\prime}(\theta)d\theta\Bigr)\Bigr).$
\end{theo}

\noindent{\bf Proof.}

 The Peano Kernel method with $n = 2,$ and $a = m,$ gives
 $$g(X) = \sum_{k = 0}^2 {(X - m)^k\over k!}g^{(k)}(m) + 
   {1\over 2!}\int_m^X (X - \theta)^2g^{\prime\prime\prime}(\theta)d\theta.$$
 Taking expected values, we have
 $$E[Y] = E[g(X)] = g(m) + {g^{\prime\prime}(m)\over 2!}\cdot\hbox{Var}[X]
          + {1\over 2!}E\Bigl(\int_m^X (X - \theta)^2
            g^{\prime\prime\prime}(\theta)d\theta\Bigr).$$
 Now, the assertion of this theorem follows from the last equation.
 
\vrule height 8pt width 4pt

 The counter examples given in [4] have this condition satisfied. The
 theorem also provides us to extend the statistical differential
 approximation beyond the third term, by requiring appropriate number of
 central moments about mean.

\begin{theo}
 If $g(\cdot)$ has derivatives of order $j \leq n,$ $g \in C^{n + 1}[I],$ where $I$ is
 some interval containing the range of $X,$ and $E\Bigl(\int_m^X(X - \theta)^n
 g^{(n +1)}(\theta)d\theta\Bigr) = 0,$ 
 then the statistical differential approximations for the expected value 
 of the transformations of the random variable $X,$ $E[g(X)]$ is exact.
 Otherwise, it will be exact up to the error term 
 $O\Bigl(E\Bigl(\int_m^X (X - \theta)^n
         g^{(n + 1)}(\theta)d\theta\Bigr)\Bigr).$
\end{theo}

 One of the short-coming of the last two theorems is that depending on the
 nature of the function $g(\cdot),$ the verification of this condition may
 be just as difficult as finding $E[g(X)]$ in some cases.

\vskip 1em
\noindent{\bf IV.~EXACTNESS OF THE VARIANCE:}
\vskip 1em

 By definition,
 $$\hbox{Var}[Y] = E\Bigl[\Bigl(g(X) - E\Bigl(g(X)\Bigr)\Bigr)^2\Bigr].$$
 Truncating the Taylor series for $g(\cdot),$ depending on the number of 
 terms of the approximating required, gives $E[Y],$ which then will be used 
 to find $\hbox{Var}\Bigl[g(X)\Bigr].$

 Suppose for an example, the expectation of $g(X),$
 $$E[Y] \cong g(m) + {1\over 2}g^{\prime\prime}(m)
              \cdot\hbox{Var}[X]$$
 is used, then the variance of $g(X)$ is
 $$\hbox{Var}[Y] = E[(Y - E[Y])^2] \cong E\Bigl\{g(X) - g(m) - 
   {1\over 2}g^{\prime\prime}(m)\cdot\hbox{Var}[X]\Bigr\}^2.$$
 Hence, the exactness of the variance formula still holds, if the function
 $g(\cdot),$ and the random variable $X$ satisfy the condition stipulated 
 in the theorems.
 
 Subject to first two terms of the expression for $E[g(X)],$ we have
 $$\hbox{Var}[Y] \cong \bigl[g^\prime(m)\bigr]^2\cdot\hbox{Var}[X].$$

 The applicable multivariate versions, involving covariance etc. given in 
 [3] can also be derived in a similar manner. The necessary steps and conditions 
 in deriving the approximate expression for $\hbox{Var}[Y]$ are similar to those
 considered in the preceding discussion.

\vskip 1em
\noindent{\bf V.~ERROR BOUNDS:}
\vskip 1em
 
 Since there are only few instances, where the approximation holds to be
 exact, we have no alternative, but to obtain some error bounds for this
 approximation. Then, it would be a question of deciding how small these
 bounds are. The following bounds are obtained, so that the accuracy of the
 approximations now entirely depend on the smallness of the error bounds,
 so desired.

 Some of the error bounds are computed for a class of functions,
 ${\cal L},$ defined by
 $${\cal L}_g \equiv \bigl\{g : g \in C^n[m, b]~ \&
            ~ |g^{(j)}(\cdot)| \leq |g(\cdot)|^{(j)}
            ~\hbox{for all}~ j \geq 1\bigr\}.$$
 This means that the derivatives are invariant under absolute value function.
 Most of the functions considered in the literatures belong to this class.
 Note that $g(x) = {1\over x^\alpha}$ for $ x > 0,$ and $\alpha > 0,$ does
 not belong to this class.

\begin{theo}
 If $g(\cdot)$ has all derivatives of order $n \geq 1$ such that 
 $|g^{(j)}(\cdot)| \leq |g(\cdot)|^{(j)},$ for all $j \geq 3,$ and 
 $E\Bigl(|X -m|^3|g^{\prime\prime\prime}(X)|\Bigr)$ exists, then the 
 statistical differential approximations for the expected value of the
 transformations of the random variable $X,$ $E[g(X)]$ is exact up to the
 error term $O\Bigl(E\bigl((X -m)^3|g^{\prime\prime\prime}(X)|\bigr)\Bigr).$
\end{theo}

\noindent{\bf Proof.}

 From the Taylor expansion about mean, we have,
\begin{eqnarray}\label{2}
 &   & g(X) - g(m) - (X - m)g^\prime(m) - 
      {(X - m)^2\over 2!}g^{\prime\prime}(m)\nonumber\\
 & = & {(X - m)^3\over 3!}g^{\prime\prime\prime}(m)
    + {(X - m)^4\over 4!}g^{(IV)}(m) + \cdots + 
      {(X - m)^n\over n!}g^{(n)}(m) + \cdots\nonumber\\
 & = & 
 (X - m)^3\Biggl\{\sum_{j = 0}^{n - 4}{(X - m)^j\over (j +3)!}(g^{\prime
    \prime\prime}(m))^{(j)} + 
    {(X -m)^{n -3}\over n!}(g^{\prime\prime\prime}(m))^{(n -3)} +
    \cdots\Biggr\}\nonumber\\
 & \leq & (X - m)^3\Biggl\{\sum_{j = 0}^{n -4}{(X -m)^j\over (j +3)!}
    |(g^{\prime\prime\prime}(m))^{(j)}| + 
    {(X -m)^{n -3}\over n!}|(g^{\prime\prime\prime}(m))^{(n -3)} +
    \cdots|\Biggr\}\nonumber\\
 & = & (X - m)^3\Biggl\{\sum_{j = 0}^{n - 4}{j!\over (j + 3)!}
     {(X - m)^j\over j!}
    |(g^{\prime\prime\prime}(m))^{(j)}| + {(n -3)!\over n!}
      {(X -m)^{n -3}\over (n -3)!}
    |(g^{\prime\prime\prime}(m))^{(n -3)}| + \cdots\Biggr\}\nonumber\\
 & \leq & (X - m)^3\Biggl\{\sum_{j = 0}^{n - 4}{(X - m)^j\over j!}
    |(g^{\prime\prime\prime}(m))^{(j)}|  + {(X -m)^{n -3}\over (n -3)!}
    |(g^{\prime\prime\prime}(m))^{(n -3)}| + \cdots\Biggr\}.
\end{eqnarray}
 Using the fact that $|g^{(j)}(m)| \leq |g(m)|^{(j)},$ we obtain,
\begin{eqnarray*}
 (\ref{2}) & \leq &
    (X - m)^3\Biggl\{\sum_{j = 0}^{n - 4}{(X - m)^j\over j!}
    (|g^{\prime\prime\prime}(m)|)^{(j)} + {(X -m)^{n -3}\over (n -3)!}
    (|g^{\prime\prime\prime}(m)|)^{(n -3)} + \cdots\Biggr\}\\
 & = & (X - m)^3|g^{\prime\prime\prime}(X)|.
\end{eqnarray*}
 Taking the expected values, we obtain the required error bound as in the
 theorem.

\vrule height 8pt width 4pt

 The following corollary follows easily.

\begin{coro}
 If $g(\cdot)$ has bounded derivatives of order $n \geq 1,$ and
 $E\Bigl(|X -m|^3e^{|X -m|}\Bigr)$ exists, then the statistical 
 differential approximations for the expected value of the
 transformations of the random variable $X,$ $E[g(X)]$ is exact up to
 the error term $O\Bigl(E\bigl((X - m)^3e^{(X -m)}\bigr)\Bigr).$
\end{coro}

\noindent{\bf Proof.}

 The proof of this corollary easily follows from the proof of the above 
 theorem.

\vrule height 8pt width 4pt

 In the following, we use the monotonicity of the norm property of random
 variables, in a certain fashion, so as to find an error bound for this
 approximation. This is in fact, so called {\sl Lyapunov Inequality.} For
 random variable $X,$
 $$\bigl\{E[|X|^s]\bigr\}^{1/s} \geq \bigl\{E[|X|^r\bigr\}^{1/r},
   ~\hbox{for all}~0 < r < s.$$
 In our case, for random variable $X - m,$ and $0 < j < n,$ we have
 $$\bigl\{E[|X - m|^n]\bigr\}^{j/n} \geq \bigl\{E[|X - m|^j\bigr\},
   ~\hbox{for all}~0 < j < n.$$

 Thus, if we assume a particular higher absolute central moment about mean
 equals zero, then all other lower absolute central moments about mean will
 be zero. This case is not much interest to us as the approximation 
 collapses to a trivial case. One of the required conditions, namely the
 condition (ii);
 $E[(X -m)^3] = 0$ implies $E[(X -m)^k] = 0,$ for all $k > 3$ 
 has no rigorous impact in this case too. However, using the monotonicity
 property, we can establish an easy result.

\begin{theo}
 If $g(\cdot)$ has all derivatives of order $n \geq 1$ such that
 $g^{(j)}(\cdot) \geq 0,$ for all $n \geq j \geq 3,$ and 
 $E[|X -m|^n]$ exists for some integer $n \geq 1,$ and let 
 $C = \bigl\{E[|X - m|^n]\bigr\}^{1/n} < \infty,$ then the statistical 
 differential approximations for the expected value of the transformations
 of the random variable $X,$ $E[g(X)]$ is exact up to the error term 
 $O\Bigl(g^{\prime\prime\prime}(C + m)\Bigr).$
\end{theo}

\noindent{\bf Proof.}

 From the Taylor expansion about mean, we have
\begin{eqnarray*}
 &   & 
 g(X) - g(m) - (X - m)g^\prime(m) + {(X -m)^2\over 2!}g^{\prime\prime}(m)\\
 & = & {(X -m)^3\over 3!}g^{\prime\prime\prime}(m) + {(X - m)^4\over 4!}
     g^{(IV)}(m) + \cdots + {(X -m)^n\over n!}g^{(n)}(z),
 ~\hbox{where}~z~\hbox{lies between}~m~\hbox{and}~b.
\end{eqnarray*}
 Evaluating for the absolute values of the right side of the expression, by
 using properties of $g^{(j)}(\cdot),$ for all $n \geq j \geq 3,$ and then 
 taking the expected values,
\begin{eqnarray*}
 &  &
   E[g(X)] - g(m) - {1\over 2}g^{\prime\prime}(m)\cdot\hbox{Var}[X]\\
 & \leq & {E|X -m|^3\over 3!}(g^{\prime\prime\prime}(m))^{(0)} +
        {E|X -m|^4\over 4!}(g^{\prime\prime\prime}(m))^{(1)} + \cdots +
        E\Biggl\{{|X -m|^n\over n!}(g^{\prime\prime\prime}(z))^{(n - 3)}\Biggr\}.
\end{eqnarray*}
 Now, using Lyapunov Inequality, we have
 $$\leq {\{E|X -m|^n\}^{3/n}\over 3!}(g^{\prime\prime\prime}(m))^{(0)} +
        {\{E|X -m|^n\}^{4/n}\over 4!}(g^{\prime\prime\prime}(m))^{(1)}
        + \cdots +
        E\Biggl\{{\{|X -m|^n\}^{n/n}\over n!}(g^{\prime\prime\prime}(z))^{(n - 3)}\Biggr\}.$$
 In summation notation, this equals to
\begin{eqnarray*}
 & = & \{E|X - m|^n\}^{3/n}\Biggl\{\sum_{j = 0}^{n - 4}{\{E|X - m|^n\}^{j/n}
    \over (j +3)!}(g^{\prime\prime\prime}(m))^{(j)} + 
    {\{E|X -m|^n\}^{(n -3)/n}\over n!}E(g^{\prime\prime\prime}(z))^{(n -3)}
    \Biggr\}\\
 & \leq & 
   \{E|X - m|^n\}^{3/n}\Biggl\{\sum_{j = 0}^{n - 4}{\{E|X - m|^n\}^{j/n}
    \over j!}(g^{\prime\prime\prime}(m))^{(j)} + 
    {\{E|X -m|^n\}^{(n -3)/n}\over (n - 3)!}
     E(g^{\prime\prime\prime}(z))^{(n -3)}\Biggr\}\\
 & = &  \{E|X - m|^n\}^{3/n}g^{\prime\prime\prime}(C + m)\\
 & = & C^3g^{\prime\prime\prime}(C + m)).
\end{eqnarray*}
 This gives the required error bound having determined an interval containing $z$
 that is independent of the range of $X.$

\vrule height 8pt width 4pt

 Further, assuming that all derivatives of $g(\cdot)$ are bounded, an corollary 
 is immediate.

\begin{coro}
 If $g(\cdot)$ has bounded derivatives of order $n \geq 1,$ and 
 $e^{\bigl(E[|X -m|^n]\bigr)^{3/n}}$ exists for some integer $n \geq 1,$
 and is finite, then the statistical differential 
 approximations for the expected value of the transformations of the random 
 variable $X,$ $E[g(X)]$ is exact up to the error term 
 $O\Bigl(\bigl\{E[|X -m|^n]\bigr\}^{3/n}
  e^{\bigl(E[|X -m|^n]\bigr)^{1/n}}\Bigr).$
\end{coro}

\noindent{\bf Proof.}

 The proof of this is essentially similar to the proof of the last theorem.

\vrule height 8pt width 4pt

 All of these results derived above can be extended for the statistical
 differential approximation truncated beyond third term under appropriate
 conditions. 

\vskip 1em
\noindent{\bf SUMMARY:}
\vskip 1em

 This topic has been presented in a course on survival methods as a prelude 
 to the other relevant chapters to follow, but failed to discuss the 
 conditions under which the exactness of this approximation to hold. The 
 lack of them, students would wonder is this approximation reasonable ?, and
 are also eager to find out the validity and accuracy of these results. A
 part of this discussion enriches rather subtle, and interesting topic, thus
 requiring inclusion in the future additions of [3]. Relevant rates of 
 convergence, error analysis and similar results for other series expansions
 can be studied, if one needs to develop this topic for further research.

\vskip 1em
\noindent{\bf ACKNOWLEDGMENTS:}
\vskip 1em

 The author wishes to thank Professor William B. Frye of Ball State University for
 many discussions over this topic, which eventually motivated to realize that there 
 are more to this interesting topic. His patience for some of the questions is very 
 appreciated. Thanks are also due to the author of [5]. His comments and loud thinking
 (in his own words) immensely helped to improve this manuscript. The author is also
 grateful to the referees and many others for their careful reading, corrections and
 helpful suggestions of the manuscript. Their comments helped enormously to improve
 this manuscript.

\vskip 1em
\noindent{\bf REFERENCES:}
\vskip 1em

\begin{enumerate}

\item Norman L. Johnson and Samuel Kotz,
      {\sl Discrete Distributions,} Houghton Mifflin, Boston (1969).

\item Regina C. Elandt-Johnson and Norman L. Johnson, 
      {\sl Survival Models and Data Analysis,} John Wiley \& Sons,
       New York (1980).

\item Dick London, {\sl Survival Models and Their estimation,}
      ACTEX Publications Inc., Winsted, Connecticut (1998).

\item William B. Frye, Method  of Statistical Differentials, 
      {\sl PRIMUS,} Volume VII, No. 3, pp. 271 - 276 (1997).

\item Kazim M. Khan, {\sl Probability with Applications,}
      MAKTABA-TUL-ILMIYA, Lahore, Pakistan (1994).

\end{enumerate}
\noindent{March 31, 2004}\\
\noindent{Department of Mathematical and Physical Sciences}\\
\noindent{Texas A\&M International University}\\
\noindent{Laredo, Texas 78041-1900}\\
\noindent{U.S.A.}
\end{document}